\documentclass[reqno]{amsart}
\usepackage[dvipsnames]{xcolor}
\usepackage{amssymb,amscd,bbm,verbatim,graphicx,mathtools}
\usepackage[mathscr]{euscript}
\usepackage[colorlinks]{hyperref}
\newtheorem{thm}{Theorem}[section]

\newtheorem{lem}[thm]{Lemma}

\theoremstyle{definition}

\newtheorem{hyp}[thm]{Hypothesis}

\numberwithin{equation}{section}

\newcommand{\supp}{\operatorname{supp}}

\newcommand{\ran}{\operatorname{ran}}

\newcommand{\diag}{\operatorname{diag}}
\newcommand{\sm}[1]{\big(\begin{smallmatrix}#1\end{smallmatrix}\big)}
\newcommand{\bb}[1]{{\mathbb{#1}}}
\newcommand{\mc}[1]{{\mathcal{#1}}}

\newcommand{\id}{\mathbbm 1}

\newcommand{\loc}{{\rm loc}}
\newcommand{\mx}{{\rm max}}
\newcommand{\mn}{{\rm min}}

\newcommand{\<}{\langle}
\renewcommand{\>}{\rangle}

\newenvironment{reduce}
 {\hbox\bgroup\scriptsize$\displaystyle}
 {$\egroup}

\begin{document}
\title[Spectral theory]{On the spectral theory for first-order systems without the unique continuation property}
\author{Kevin Campbell, Minh Nguyen, and Rudi Weikard}
\address{Department of Mathematics, University of Alabama at Birmingham, Birmingham, AL 35226-1170, USA}
\email{campke@uab.edu, minhnt@uab.edu, weikard@uab.edu}
\date{3. October 2019}
\thanks{This is an Accepted Manuscript of an article to be published by Taylor \& Francis in Linear and Multilinear Algebra, available online at \texttt{https://doi.org/10.1080/03081087.2019.1671303}.}
\thanks{\copyright 2019. This manuscript version is made available under the CC-BY-NC-ND 4.0 license \texttt{http://creativecommons.org/licenses/by-nc-nd/4.0/}}

\begin{abstract}
We consider the differential equation $Ju'+qu=wf$ on the real interval $(a,b)$ when $J$ is a constant, invertible skew-Hermitian matrix and $q$ and $w$ are matrices whose entries are distributions of order zero with $q$ Hermitian and $w$ non-negative.
In this situation it may happen that there is no existence and uniqueness theorem for balanced solutions of a given initial value problem.
We describe the set of solutions the equation does have and establish that the adjoint of the minimal operator is still the maximal operator, even though unique continuation of balanced solutions fails.
\end{abstract}
\maketitle

\section{Introduction}
Ghatasheh and Weikard \cite{GW18-1} investigated the spectral theory for the first-order system
$$Ju'+qu=wf$$
of differential equations on the real interval $(a,b)$ assuming that $J$ is a constant, invertible skew-Hermitian matrix and $q$ and $w$ are matrices whose entries are distributions of order zero\footnote{Recall that
distributions of order $0$ are distributional derivatives of functions of locally bounded variation and hence may be thought of, on compact subintervals of $(a,b)$, as measures. For simplicity we might use the word measure instead of distribution of order $0$ below.} with $q$ Hermitian and $w$ non-negative.
Crucially, \cite{GW18-1} requires that initial value problems for this equation have unique balanced\footnote{The concept of balanced solutions is defined below and the rationale of using it is explained in \cite{GW18-1}.} solutions.
Indeed, unique continuation of a solution across a point where $Q$, the anti-derivative of $q$, has a discontinuity may fail.
With the aid of linear algebra we were able to overcome this obstacle, describe the set of solutions of the differential equation, and establish that the adjoint of the minimal operator is still the maximal operator, even if unique continuation of balanced solutions fails.

The relationship between minimal and maximal operators is a cornerstone for the spectral theory for differential equations.
Of course, this topic is very well studied when the coefficients are locally integrable functions but in the case of measure coefficients much less is known.
The first to consider an equation with a measure coefficient was Krein \cite{MR0054078} in 1952 when he modeled a vibrating string.
Also motivated by physical applications were Gesztesy and Holden \cite{MR914699} in 1987 who described Schr\"odinger equations with point interactions, specifically $\delta'$-interactions.
In 1999 Savchuk and Shkalikov \cite{MR1756602} treated Schr\"odinger equations with potentials in the Sobolev space $W^{-1,2}_{\rm\loc}$, a paper which spurred many further developments.
With the help of quasi-derivatives Eckhardt et al. \cite{MR3046408} showed in 2013 that such equations can be cast as first-order $2\times2$-systems with locally integrable coefficients.
Eckhardt and Teschl \cite{MR3095152} considered a system where the coefficients are measures, viz. $Ju'+qu=wf$ where $J=\sm{0&-1\\1&0}$, $q=\sm{\chi&0\\0&-\varsigma}$, and $w=\sm{\rho&0\\0&0}$.
Their approach covers both the Krein string ($\chi=0$ and $\varsigma=1$) as well as the $\delta'$-interaction ($\chi=0$, $\varsigma=1+\beta\delta_0$, and $\rho=1$ in the simplest case).
Crucially, they require that the support of the discrete part of $\varsigma$ does not intersect the corresponding sets for $\chi$ or $\rho$, a condition which guarantees unique continuation.
Both \cite{MR3046408} and \cite{MR3095152} are also excellent sources for a more thorough history of the subject.

Let us add a few words about notation.
$\mc D^{\prime0}((a,b))$ is the space of distributions of order $0$, i.e., the space of distributional derivatives of functions of locally bounded variation.
Any function $u$ of locally bounded variation has left- and right-hand limits denoted by $u^-$ and $u^+$, respectively.
Also, $u$ is called balanced if $u=u^\#=(u^++u^-)/2$.
We use $\id$ to denote an identity matrix of appropriate size and superscripts ${}^\top$ and ${}^*$ indicate transposition and adjoint, respectively.
The orthogonal complement of a subspace $S$ of a Hilbert space $H$ is denoted by $H\ominus S$ or by $S^\perp$.
For $c_1, ..., c_N\in\bb C^n$ we abbreviate the vector $(c_1^\top, ..., c_N^\top)^\top\in\bb C^{nN}$ by $(c_1, ..., c_N)^\diamond$.
Finally we note that, generally, our differential equations are represented by linear relations rather than linear operators.
Consequently we work with graphs of such relations (even when they are operators).

\section{Obtaining solutions}\label{s2}
We begin by describing the set of solutions of the first-order system
$$Ju'+qu=wf$$
on the interval $(a,b)$ assuming that the coefficients satisfy the following hypothesis.
\begin{hyp}\label{H:1}
$J$ is a constant, invertible and skew-Hermitian $n\times n$-matrix.
Both $q$ and $w$ are in $\mc D^{\prime0}((a,b))^{n\times n}$, $w$ is non-negative, and $q$ Hermitian.
\end{hyp}

Associated with a non-negative distribution $w\in \mc D^{\prime0}((a,b))^{n\times n}$ is a Hilbert space $L^2(w)$ with inner product $\<u,v\>=\int u^*wv$ (recall that positive distributions are positive measures).
Its elements are equivalence classes of functions $[f]$ satisfying $\|f\|^2=\int f^*wf<\infty$ and, as usual, two functions $f$ and $g$ are equivalent, if $\|f-g\|=0$.

We denote the left-continuous anti-derivatives of $q$ and $w$ by $Q$ and $W$, respectively.
$\Delta_q(x)=Q^+(x)-Q^-(x)$ stands for the jump of $Q$ at a point $x$.
Similarly, $\Delta_w(x)=W^+(x)-W^-(x)$.

Suppose $\xi_0\in(\xi_1,\xi_2)\subset(a,b)$.
When $f\in L^2(w)$ (as we shall henceforth assume) it was shown in \cite{GW18-1} that the initial value problem $Ju'+qu=wf$, $u(\xi_0)=u_0\in\bb C^n$ has a unique balanced solution of locally bounded variation in $(\xi_1,\xi_2)$ provided that the matrices $B_\pm(x)=J\pm\Delta_q(x)/2$ are invertible for all $x\in(\xi_1,\xi_2)$.
At a point $x$ of discontinuity of $Q$ or $W$, the differential equation requires that
$$J(u^+(x)-u^-(x))+\Delta_{q}(x)u(x)=\Delta_w(x)f(x)$$
where, $u$ being balanced, $u(x)=(u^+(x)+u^-(x))/2$.
This is equivalent to
\begin{equation}\label{e:2.1}
B_+(x)u^+(x)-B_-(x) u^-(x)=\Delta_w(x)f(x).
\end{equation}
From this it is obvious that we may not be able to continue a solution across $x$ from left to right (or from right to left), if $B_+(x)$ (or $B_-(x)$) fails to be invertible.
In our particular situation, where $J^*=-J$ and $\Delta_q(x)^*=\Delta_q(x)$, we have $B_-(x)=-B_+(x)^*$ and hence that $B_-(x)$ is invertible if and only if $B_+(x)$ is.

On account of the fact that $Q$ is locally of bounded variation, it is clear that the set of points $x$ where $B_\pm(x)$ are not invertible is discrete and finite on compact subintervals of $(a,b)$ even though the set of all jumps of $Q$ may be dense.

Let us now fix an interval $[\xi_1,\xi_2]\subset (a,b)$ assuming that the points in $(\xi_1,\xi_2)$ where $B_\pm$ are not invertible are among the points $x_1< ... <x_N$.
Normally one would choose these to be precisely the points where $B_\pm$ are not invertible but it is advantageous to avoid the case $N=1$.
For convenience let us also set $x_0=\xi_1$ and $x_{N+1}=\xi_2$.
As mentioned above we do have unique solutions of initial value problems and, indeed, a variation of constants formula in any of the intervals $(x_j,x_{j+1})$.
These solutions have limits at the endpoints of the interval and we may even use, for instance, the left endpoint to pose an initial condition.
Therefore the general solution of $Ju'+qu=wf$ in $(x_j,x_{j+1})$ is represented by
\begin{equation}\label{e:2.1.a}
u^-(x)=U_j^-(x)(c_j+J^{-1} \int_{(x_j,x)} U_j^*wf)
\end{equation}
where $c_j$ is an arbitrary element of $\bb C^n$ and $U_j$ a balanced fundamental matrix for $Ju'+qu=0$ in $(x_j,x_{j+1})$ which we may choose so that $\lim_{x\downarrow x_j}U_j(x)=\id$.
We then define $U_j(x_{j+1})=\lim_{x\uparrow x_{j+1}} U_j^-(x)$.
For $u$ to be a solution of $Ju'+qu=wf$ on $(x_0,x_{N+1})$ we need $u$ to be determined by \eqref{e:2.1.a} (for appropriate choices of the $c_j$) in the respective intervals.
Moreover, according to equation \eqref{e:2.1}, $u$ must satisfy
\begin{equation}\label{e:2.2}
B_+(x_j)u^+(x_j)-B_-(x_j) u^-(x_j)=\Delta_w(x_j)f(x_j) \quad \text{for $j=1,..., N$}.
\end{equation}
Note that $u^+(x_j)=c_j$ and $u^-(x_j)=U_{j-1}(x_j)(c_{j-1}+J^{-1} I_{j-1}(f))$ where $I_{j-1}(f)=\int_{(x_{j-1},x_j)} U_{j-1}^*wf$.
Thus we may rewrite equation \eqref{e:2.2} as
$$(-B_-(x_j)U_{j-1}(x_j), B_+(x_j))\sm{c_{j-1}\\ c_j}= \Delta_w(x_j)f(x_j)+B_-(x_j)U_{j-1}(x_j)J^{-1} I_{j-1}(f).$$

At this point it appears helpful to introduce the following notation.
Let
\begin{eqnarray*}
  \mc B &=& \diag(B_+(x_1),..., B_+(x_N)), \\
  \mc U &=& \diag(U_0(x_1),..., U_{N-1}(x_N)),\\
  \mc J &=& \diag(J, ..., J),
\end{eqnarray*}
and $E_\top=(0,\id)$ and $E_\bot=(\id,0)$, two $nN\times n(N+1)$-matrices which, respectively, strip the first and the last $n$ coordinates off a vector.
Then we have
\begin{multline*}
B = \mc B^*\mc U E_\bot+\mc B E_\top \\
=\begin{reduce}
\begin{pmatrix}
-B_-(x_1)U_0(x_1)&B_+(x_1)&0& \cdots &0\\
0&-B_-(x_2)U_1(x_2)&B_+(x_2)& \cdots &0\\
\vdots&\vdots&\vdots&\ddots&\vdots\\
0&\cdots&0&-B_-(x_N)U_{N-1}(x_N)&B_+(x_N)
\end{pmatrix}
\end{reduce}.
\end{multline*}
If we now introduce the abbreviations
\begin{eqnarray*}
\tilde u &=& (c_0, ...,c_N)^\diamond,\\
\mc R(f) &=& (\Delta_w(x_1)f(x_1), ..., \Delta_w(x_N)f(x_N))^\diamond,\\
\mc I(f) &=& (I_0(f), ..., I_{N-1}(f))^\diamond,
\end{eqnarray*}
and, for later purposes,
$$\tilde{\mc I}(f)=(0,...,0,I_N(f))^\diamond\in\bb C^{nN},$$
equations \eqref{e:2.2} may be written as
\begin{equation}\label{e:2.3}
B\tilde u = \mc R(f)-\mc B^*\mc U\mc J^{-1} \mc I(f).
\end{equation}

We have proved the following result.
\begin{thm}\label{mansol}
If $u$ is any solution of $Ju'+qu=wf$ on $(x_0,x_{N+1})$ then $\tilde u=(u^+(x_0), ..., u^+(x_N))^\diamond$ is a solution of equation \eqref{e:2.3}.
Conversely, a solution $\tilde u=(c_0,...,c_N)^\diamond$ of equation \eqref{e:2.3} provides a solution of $Ju'+qu=wf$ on $(x_0,x_{N+1})$ given by \eqref{e:2.1.a} for $j=0,..., N$.
\end{thm}

It is clear that the rank of $B$ is no larger than $nN$ and hence the kernel of $B$ has at least dimension $n$.
Note, however, that it is possible for the dimension of the kernel of $B$ to be larger than $n$, i.e., to have more than $n$ independent solutions of the homogeneous differential equation $Ju'+qu=0$.

When we consider balanced solutions of $Ju'+qu=0$, the relationship between $\tilde u$ and the vector of values of $u$ at the points $x_1$, ..., $x_N$, i.e., the vector $\hat u=(u(x_1),..., u(x_N))^\diamond$ is given by $\hat u=C\tilde u$ where
$$C = \frac12 (\mc U E_\bot+E_\top)
=\frac12\begin{reduce}
\begin{pmatrix}
U_0(x_1)&\id&0& \cdots &0\\
0&U_1(x_2)&\id& \cdots &0\\
\vdots&\vdots&\vdots&\ddots&\vdots\\
0&\cdots&0&U_{N-1}(x_N)&\id
\end{pmatrix}
\end{reduce}.$$
We also introduce the matrices $B_m$ and $C_m$ which are obtained from $B$ and $C$ respectively by removing the first and last $n$ columns.
Earlier we chose, without loss of generality, $N\geq2$ to avoid the case when $B_m$ and $C_m$ have no columns.

We have the following relationship between $B$ and $C$.
\begin{lem}\label{cbbc}
$C^*B-B^*C=\diag(-J,0,...,0,J)$.
In particular, $C_m^*B-B_m^*C=0$.
\end{lem}
\begin{proof}
Since $\mc B-\mc B^*=2\mc J$ we get
$$C^*B-B^*C=(E_\top^*\mc JE_\top-E_\bot^*\mc U^*\mc J \mc UE_\bot).$$
It was shown in \cite{GW18-1} that $u^{-*}Jv^-$ is constant on any interval on which $B_\pm$ are everywhere invertible when $u$ and $v$ are solutions of $Ju'+qu=0$.
In particular, $U_j(x_{j+1})^*JU_j(x_{j+1})=\lim_{x\uparrow x_{j+1}}U_j^{-}(x)^*JU_j^-(x)=J$.
This implies that $\mc U^*\mc J \mc U=\mc J$ and hence the claim.
\end{proof}

\begin{lem}\label{findtu}
Suppose $\hat u\in\ker B_m^*$. Then there exists a unique vector $\tilde u$ such that $B\tilde u=0$ and $C\tilde u=\hat u$.
Moreover, if $\hat u\in\ker B^*\subset\ker B_m^*$, then the first and the last $n$ components of $\tilde u$ are equal to $0$.
\end{lem}
\begin{proof}
If a solution $\tilde u$ indeed exists, it must satisfy $\mc U E_\bot\tilde u=2\hat u-E_\top\tilde u$ and hence $0=2\mc B^*\hat u +(\mc B-\mc B^*)E_\top\tilde u$.
This implies $E_\top\tilde u=-\mc J^{-1}\mc B^*\hat u$.
Similarly, using $E_\top\tilde u=2\hat u-\mc UE_\bot\tilde u$, we get $E_\bot\tilde u=\mc J^{-1}\mc U^*\mc B\hat u$.
Thus a solution is unique.

To prove existence note that $\hat u=(\hat u_1, ...,\hat u_N)^\diamond\in\ker B_m^*$ implies
$$B_+(x_k)^*\hat u_k+U_{k}(x_{k+1})^*B_+(x_{k+1})\hat u_{k+1}=0$$
for $k=1,...,N-1$.
Hence the assignments $E_\top\tilde u=-\mc J^{-1}\mc B^*\hat u$ and $E_\bot\tilde u=\mc J^{-1}\mc U^*\mc B\hat u$ define $\tilde u$ unambiguously.
Also $\tilde u$ satisfies $B\tilde u=0$ and $C\tilde u=\hat u$.

For the last claim notice that $(C^*B-B^*C)\tilde u=0$.
If $\tilde u=(c_0,..., c_N)^\diamond$, Lemma~\ref{cbbc} gives $-Jc_0=Jc_N=0$ and hence $c_0=c_N=0$ as claimed.
\end{proof}

\begin{thm}\label{T:2.5}
If $\hat u\in\ker B_m^*$, then $Ju'+qu=0$ has a unique solution $u$ on the interval $(x_0,x_{N+1})$ such that $\tilde u=(u^+(x_0),...,u^+(x_N))^\diamond$ satisfies $B\tilde u=0$ and $C\tilde u=(u(x_1),..., u(x_N))^\diamond=\hat u$.
If $\hat u\in\ker B^*\subset \ker B_m^*$ then, additionally, $u^+(x_0)=u^-(x_{N+1})=0$ so that $\supp u\in[x_1,x_N]$.
\end{thm}

\begin{proof}
This is an immediate consequence of Theorem \ref{mansol} and Lem\-ma~\ref{findtu}.
\end{proof}

\begin{lem}\label{L:2.5}
Suppose $f\in L^2(w)$ and $\hat u\in\ker B_m^*$.
Then there is a function $u$ satisfying $Ju'+qu=0$ on $(x_0,x_{N+1})$, $(u(x_1), ..., u(x_N))^\diamond=\hat u$, and
$\hat u^*\mc F(f)=\int_{(x_0,x_{N+1})} u^*wf$,
where
$$\mc F(f)=\mc R(f)-\mc B^*\mc U\mc J^{-1} \mc I(f)+\mc B\mc J^{-1}\tilde{\mc I}(f).$$
\end{lem}

\begin{proof}
Let $u$ be the function furnished by Theorem \ref{T:2.5} and let $\tilde u$ be the vector $(u^+(x_0), ..., u^+(x_N))^\diamond$.
Then $\mc B E_\top\tilde u=-\mc B^*\mc UE_\bot\tilde u$ since $B\tilde u=0$.
We also have $2C\tilde u=\mc U E_\bot\tilde u+E_\top\tilde u$.
Using $\mc B-\mc B^*=2\mc J$ gives us that $\mc B C\tilde u=\mc J\mc U E_\bot\tilde u$ and $\mc B^* C\tilde u=-\mc J E_\top\tilde u$
or, taking adjoints,
$$\tilde u^*C^*\mc B^*=-\tilde u^*E_\bot^*\mc U^*\mc J\quad\text{and}\quad \tilde u^*C^*\mc B=\tilde u^*E_\top^*\mc J.$$
These and $\mc U^*\mc J\mc U=\mc J$ imply
\begin{multline*}
\hat u^*\mc F(f) =\hat u^* \mc R(f)-\tilde u^*C^*\mc B^*\mc U \mc J^{-1}\mc I(f)+\tilde u^*C^*\mc B\mc J^{-1}\tilde{\mc I}(f)\\
=\hat u^* \mc R(f)+\tilde u^* E_\bot^*\mc I(f)+\tilde u^*E_\top^*\tilde{\mc I}(f) \\
=\hat u^* \mc R(f)+\tilde u^*(I_0(f),...,I_N(f))^\diamond =\int_{(x_0,x_{N+1})} u^*wf
\end{multline*}
using that $u(x)=U_j(x)u^+(x_j)$ for $x\in(x_j,x_{j+1})$.
\end{proof}

\section{Minimal and maximal relations}\label{s3}
Our differential equation $Ju'+qu=wf$ gives rise to the following two linear relations.
$T_\mx$ is the set of all pairs $([u],[f])\in L^2(w)\times L^2(w)$ for which there are representatives $u\in [u]$ and $f\in[f]$ such that $Ju'+qu=wf$ (in particular, $u$ is a balanced function of locally bounded variation).
$T_\mn$ is the set of those elements in $T_\mx$ for which the solution $u\in[u]$ may be chosen with compact support.

Recall that the adjoint of a linear relation $S\subset L^2(w)\times L^2(w)$ is defined to be
$$S^*=\{(v,g)\in L^2(w)\times L^2(w):\forall (u,f)\in S: \<g,u\>=\<v,f\>\}.$$
Our main result in this paper is the following theorem.
\begin{thm}\label{T:main}
$T_\mn^*=T_\mx$.
\end{thm}

Our proof requires a little preparation with which we begin.
Suppose $[\xi_1,\xi_2]\subset (a,b)$ and consider the relation $\breve T_\mx$ associated with $J$, $q$, and $w$ but restricted to the interval $(\xi_1,\xi_2)$.
Of course, solutions of $Ju'+qu=wf$ have limits at $\xi_1$ and $\xi_2$.
We denote the restriction of $w$ to $(\xi_1,\xi_2)$ by $\breve w$ and set $T_0=\{([u],[f])\in\breve T_\mx:u^+(\xi_1)=u^-(\xi_2)=0\}$ and $K_0=\ker\breve T_\mx$.
\begin{lem}\label{L:3.1}
$\ran T_0=L^2(\breve w)\ominus K_0$.
\end{lem}

\begin{proof}
Let $[f]\in \ran T_0$ and $[r]\in K_0$.
Then $Ju'+qu=\breve wf$ for some $u$ which vanishes at $\xi_1$ and $\xi_2$ and $r$ (chosen appropriately in $[r]$) satisfies $Jr'+qr=0$.
Integration by parts shows then
$$\int f^*\breve wr=\int u^*(Jr'+qr)=0.$$
Hence $\ran T_0\subset L^2(\breve w)\ominus K_0$.

Conversely, suppose $[f]\in L^2(\breve w)\ominus K_0$.
We want to show the existence of a balanced function $u$ of bounded variation defined on $(\xi_1,\xi_2)$, vanishing at the endpoints, and satisfying $Ju'+qu=\breve wf$.
Using the notation established in Section~\ref{s2} and, in particular, Theorem \ref{mansol} we have to show the existence of a solution $\tilde u=(\gamma_0, ..., \gamma _N)^\diamond$ of equation \eqref{e:2.3} satisfying $\gamma_0=0$ and $\gamma_N=-J^{-1} I_N(f)$ (so that $u^+(\xi_1)=u^-(\xi_2)=0$).
Thus we need to find $\tilde u_0=(\gamma_1,..., \gamma_{N-1})^\diamond$ such that $B_m\tilde u_0=\mc F(f)$ where, as in Lemma \ref{L:2.5},
$$\mc F(f)=\mc R(f)-\mc B^*\mc U\mc J^{-1} \mc I(f)+\mc B\mc J^{-1}\tilde{\mc I}(f).$$
This system has a solution precisely when $\mc F(f)$ is in $\ran B_m=(\ker B_m^*)^\perp$.

Hence suppose $\hat r\in\ker B_m^*$.
The function $r$ associated with $\hat r$ according to Theorem \ref{T:2.5} is a representative of an element in $K_0$ so that $\int r^* \breve w f=0$.
But Lemma \ref{L:2.5} shows that $\hat r^*\mc F(f)=\int r^* \breve w f$ guaranteeing the existence of $u$.
\end{proof}

\begin{lem}\label{L:3.2}
If $g\in\ran T_\mn^*$, then the differential equation $Ju'+qu=wg$ has at least one solution on $(a,b)$.
\end{lem}

\begin{proof}
Let $\tau_n$, $n\in \bb Z$, be an enumeration of points in $(a,b)$ which include all points where the matrices $J\pm\Delta_q(x)/2$ are not invertible.
The labeling is such that $\tau_n<\tau_{n+1}$ and we may arrange things so that $a$ and $b$ are accumulation points, and the only ones, of the sequence $\tau_n$.

According to Theorem \ref{mansol} there is a balanced solution $v_j$ of $Ju'+qu=wg$ on $(\xi_1,\xi_2)=(\tau_{-j},\tau_j)$ (at least when $j>1$) provided $B\tilde v_j=G$ where
$$G=\mc R(g)-\mc B^*\mc U\mc J^{-1} \mc I(g)=\mc F(g)-\mc B\mc J^{-1}(0,...,0,I_N(g))^\diamond.$$
This, in turn, happens if and only if $G\in\ran B=(\ker B^*)^\perp$ which we show next.

For any $\hat r\in\ker B^*$ Theorem \ref{T:2.5} and Lemma \ref{L:2.5} show the existence of a solution $r$ of $Ju'+qu=0$ on $(\xi_1,\xi_2)$ such that $\tilde r=(r^+(x_0),...,r^+(x_N))^\diamond$ satisfies $B\tilde r=0$, $r^+(x_0)=r^+(x_N)=0$, and $\hat r^*\mc F(g)=\int_{(\xi_1,\xi_2)} r^* wg$ (here $N=2j-1$ and $x_\ell=\tau_{-j+\ell}$).
In fact, since $r$ vanishes near $x_0$ and $x_{N+1}$, we may extend it by $0$ to obtain a solution of $Ju'+qu=0$ on all of $(a,b)$.
Thus $\<r,g\>=\hat r^*\mc F(g)$ and it follows, as in the proof of Lemma \ref{L:2.5}, that $\hat r^*G=\<r,g\>-\tilde r^*E_\top^*(0,...,0,I_N(g))^\diamond$.
Since $([r],0)\in T_\mn$ we have $\<r,g\>=\<0,v\>=0$ and since $r^+(x_N)=0$ we also have $\tilde r^*E_\top^*(0,...,0,I_N(g))^\diamond=0$.
Thus $G\in\ran B$ and this guarantees the existence of $v_j$.

Now define, for any $j\geq k\geq2$, the set $A_{k,j}$ to be the collection of restrictions to $(\tau_{-k},\tau_k)$ of solutions of $Ju'+qu=wg$ on $(\tau_{-j},\tau_j)$.
According to the above the $A_{k,j}$ are non-empty and nested in the sense that $A_{k,j+1}\subset A_{k,j}$.
Each $A_{k,j}$ is an affine subspace of, say, the space of all functions defined on $(\tau_{-k},\tau_k)$ and their dimensions form, in $j$, a non-increasing sequence of non-negative integers which must eventually be constant (possibly zero).
Hence, for a sufficiently large $m$, we have that $B_k=\bigcap_{j\geq k} A_{k,j}=A_{k,m}$.
We now define inductively a sequence of functions $v_k$ whose pointwise limit is a solution of $Ju'+qu=wg$ on $(a,b)$.
For $v_2$ we choose any element of $B_2$.
Then suppose we had constructed a sequence $(v_2, ..., v_k)$ such that $v_j\in B_j$ and $v_{j-1}=v_j|_{(\tau_{-j+1},\tau_{j-1})}$.
Note that the elements of $B_k$ are restrictions of elements in $B_{k+1}$ to $(\tau_{-k},\tau_k)$ (and vice versa).
Thus we may choose for $v_{k+1}$ an element of $B_{k+1}$ which extends $v_k$ and this completes our definition of the sequence $v_k$, except that we extend each of its elements arbitrarily to $(a,b)$.
Now $v$, the pointwise limit of the $v_k$, is the desired solution of $Ju'+qu=wg$ on $(a,b)$.
\end{proof}

\begin{proof}[Proof of Theorem \ref{T:main}]
If $([v],[g])$ and $([u],[f])$ are in $T_\mx$ Lagrange's identity (cf. \cite{GW18-1}) states that
\begin{equation}\label{Lagrange}
\<v,f\>-\<g,u\>=(v^*Ju)^-(b)-(v^*Ju)^+(a).
\end{equation}
Therefore, if $([u],[f])\in T_\mn$, so that $u$ has compact support, we get $\<v,f\>=\<g,u\>$ which proves that $T_\mx\subset T_\mn^*$.

To prove $T_\mn^*\subset T_\mx$ assume that $([v],[g])\in T_\mn^*$ and let $v_0$ be a solution of $Ju'+qu=wg$ on $(a,b)$ as constructed by Lemma \ref{L:3.2}.
Next we will employ Lemma \ref{L:3.1}.
We consider an interval $[\xi_1,\xi_2]\subset (a,b)$ and define $\breve w$, $T_0$ and $K_0$ as we did there.
Given $([u],[f])\in T_0$ extend both $u$ and $f$ by $0$ to all of $(a,b)$ (denoting the extensions also by $u$ and $f$).
We then have $Ju'+qu=wf$ so that $([u],[f])\in T_\mn$ and $\<f,v\>=\<u,g\>$.
To establish a relationship between $v$ and $v_0$ we apply integration by parts to obtain
\begin{multline*}
\int_{(\xi_1,\xi_2)} f^*\breve w v=\int_{(a,b)} u^*wg=\int_{(a,b)} u^*(Jv_0'+qv_0)\\
 =\int_{(a,b)} (Ju'+qu)^*v_0=\int_{(\xi_1,\xi_2)} f^*\breve w v_0.
\end{multline*}
Thus $\int f^*\breve w(v-v_0)=0$ so that, by Lemma \ref{L:3.1}, $[v-v_0]\in K_0$ showing that $[v]$ has a representative $v=v_0+k_0$ where $Jk_0'+qk_0=0$ and hence $Jv'+qv=wg$ on $(\xi_1,\xi_2)$.
We obtain a solution on all of $(a,b)$ in a similar way as we did in the proof of Lemma \ref{L:3.2}.
We only have to modify the definition of the sets $A_{k,j}$ to specify that the solutions $u$ considered are locally representatives of $v$, i.e, that $\int_{(\tau_{-j},\tau_j)} (u-v)^*w(u-v)=0$.
\end{proof}

\section*{Acknowledgements}
We thank the anonymous referee for most helpful comments and questions.


\begin{thebibliography}{1}

\bibitem{MR3046408}
Jonathan Eckhardt, Fritz Gesztesy, Roger Nichols, and Gerald Teschl.
\newblock Weyl-{T}itchmarsh theory for {S}turm-{L}iouville operators with
  distributional potentials.
\newblock {\em Opuscula Math.}, 33(3):467--563, 2013.

\bibitem{MR3095152}
Jonathan Eckhardt and Gerald Teschl.
\newblock Sturm-{L}iouville operators with measure-valued coefficients.
\newblock {\em J. Anal. Math.}, 120:151--224, 2013.

\bibitem{MR914699}
F.~Gesztesy and H.~Holden.
\newblock A new class of solvable models in quantum mechanics describing point
  interactions on the line.
\newblock {\em J. Phys. A}, 20(15):5157--5177, 1987.

\bibitem{GW18-1}
Ahmed Ghatasheh and Rudi Weikard.
\newblock Spectral theory for systems of ordinary differential equations with
  distributional coefficients.
\newblock {\em Submitted, 46 pages}, 2018.
\newblock{\texttt{arXiv:1807.09653v2}}

\bibitem{MR0054078}
M.~G. Kre{\u\i}n.
\newblock On a generalization of investigations of {S}tieltjes.
\newblock {\em Doklady Akad. Nauk SSSR (N.S.)}, 87:881--884, 1952.

\bibitem{MR1756602}
A.~M. Savchuk and A.~A. Shkalikov.
\newblock Sturm-{L}iouville operators with singular potentials.
\newblock {\em Mathematical Notes}, 66(6):741--753, 1999.
\newblock Translated from Mat. Zametki, Vol. 66, pp. 897--912 (1999).

\end{thebibliography}

\end{document}